\documentclass{amsart}
\usepackage{amscd,amssymb}
\usepackage[all]{xy}

\newcommand{\gr}{{\operatorname{gr}\nolimits}}

\newcommand{\kar}{\operatorname{char}\nolimits}
\newcommand{\Hom}{\operatorname{Hom}\nolimits}

\renewcommand{\Im}{\operatorname{Im}\nolimits}
\newcommand{\Ker}{\operatorname{Ker}\nolimits}

\newcommand{\rrad}{\mathfrak{r}}

\newcommand{\Span}{\operatorname{span}\nolimits}

\newcommand{\Ext}{\operatorname{Ext}\nolimits}

\newcommand{\op}{{\operatorname{op}\nolimits}}

\newcommand{\id}{{\operatorname{id}\nolimits}}

\newcommand{\rad}{{\operatorname{rad}\nolimits}}

\newcommand{\G}{\Gamma}
\renewcommand{\L}{\Lambda}
\newcommand{\Z}{{\mathbb Z}}

\newcommand{\E}{{\mathcal E}}

\newcommand{\extto}{\xrightarrow}

\newcommand{\HH}{\operatorname{HH}\nolimits}

\newtheorem{lem}{Lemma}[section]
\newtheorem{prop}[lem]{Proposition}

\newtheorem*{Theorem}{Theorem}
\theoremstyle{definition}

\begin{document}

\title{Finite Hochschild cohomology without finite global dimension}

\author[Buchweitz]{Ragnar-Olaf Buchweitz}
\address{Ragnar-Olaf Buchweitz\\ 
Department of Mathematics\\
University of Toronto\\ 
Toronto, ON Canada M5S 3G3\\ 
Canada}
\email{ragnar@math.toronto.edu}
\author[Green]{Edward L. Green}
\address{Edward L. Green\\ Department of Mathematics\\
Virginia Tech\\ Blacksburg, VA 24061\\ USA}
\email{green@math.vt.edu}
\author[Madsen]{Dag Madsen}
\address{Dag Madsen\\Institutt for matematiske fag\\
NTNU\\ N--7491 Trondheim\\ Norway}
\email{dagma@math.ntnu.no}
\author[Solberg]{\O yvind
Solberg$^\dagger$}\thanks{$^\dagger$Supported by the Research council
of Norway and NTNU}
\address{\O yvind Solberg\\Institutt for matematiske fag\\
NTNU\\ N--7491 Trondheim\\ Norway}
\email{oyvinso@math.ntnu.no}
\date{\today}
\maketitle

\section*{Introduction}
In \cite{H2}, Dieter Happel asked the following question: If the
Hochschild cohomology groups $\HH^n(\G)$ of a finite dimensional
algebra $\G$ over a field $k$ vanish for all sufficiently large $n$,
is the global dimension of $\G$ finite? We give a negative answer to
this question.  Indeed, consider the finite dimensional algebras
$\L=\L_q=k\langle x,y\rangle/(x^2,xy+qyx,y^2)$ for some field $k$ and
$q$ in $k$. These algebras, intensely studied already by Rainer Schulz
\cite{S} to exhibit other pathologies, are all four dimensional, of
infinite global dimension, and, for $q\neq 0$, selfinjective.  If $q$
is not a root of unity in $k$ and non-zero, then we show that the
total Hochschild cohomology of $\L_q$ is of dimension $5$, thus
answering negatively Happel's question as promised.  Moreover, if $q$
is a root of unity, another interesting phenomenon occurs: the
Hochschild cohomology vanishes on large, but finite intervals, whose
lengths depend on the multiplicative order of $q$.

None of these phenomena can occur for \emph{commutative} algebras 
according to the companion paper \cite{AI}. Indeed, there it is shown,
as a consequence of a more general result, that a commutative finite
dimensional algebra $R$ is of finite global dimension as soon as
Hochschild cohomology vanishes just in one even and one odd degree,
thus giving a strong affirmation of Happel's question in the
commutative case. Let us point out that the family $\{\L_q\}$ contains a
few (graded) commutative specimens. For $q=-1$, the algebra
$\Lambda_{-1}$ is the commutative complete intersection
$k[x,y]/(x^2,y^2)$ and for $q=1$, the algebra $\Lambda_1$ is the
(graded commutative) exterior algebra in two variables, isomorphic to
the group ring $k[\Z_2\times\Z_2]$ in characteristic two. For $q=0$,
that is $\L_0=k\langle x,y\rangle/(x^2,xy,y^2)$, we discuss the
structure of its Hochschild cohomology in section \ref{section:q=0}.

The Hochschild cohomology ring $\HH^*(\L)$ has a ring structure given
by the Yoneda product. If $R\to A$ and $S\to A$ are two ring
homomorphisms, we denote the fibre product (pullback) by $R\times_A
S$. 

In this paper we prove the following.
\begin{Theorem}
\begin{enumerate}
\item[(a)] If $q$ is not a root of unity and non-zero, then
$\HH^*(\L)$ is a $5$-dimensional algebra over the field $k$ with
$\sum_{i\geq 0} \dim_k \HH^i(\L) t^i =2+2t+t^2$. Moreover, $\HH^*(\L)$
is isomorphic as a ring to the fibre product
\[k[z]/(z^2)\times_k \wedge^*(u_1,u_2)\]
where $z$ is in degree $0$, the elements $u_1$ and $u_2$ are in degree
$1$ and $\wedge^*(u_1,u_2)$ is the exterior algebra on the generators
$u_1$ and $u_2$.
\item[(b)] If $q$ is a root of unity, then $\HH^*(\L)$ is an infinite
dimensional algebra over the field $k$ and finitely generated as a
module over a polynomial subring generated by two elements.
\item[(c)] \sloppy For $q=0$, then $\HH^*(\L)$ is an infinitely generated
algebra over $k$ with $\sum_{i\geq 0} \dim_k \HH^i(\L) t^i =
\frac{1+t^3}{(1-t)^2}$. 
\end{enumerate}
\end{Theorem}
The Hilbert series for the algebras in case (b) are found in sections
\ref{section:2.1} through \ref{section:2.5}. 

If $k\subseteq k'$ is a field extension, then
$\HH^*(k'\otimes_k\L)\simeq k'\otimes_k\HH^*(\L)$ as graded
$k'$-algebras; see, for example, \cite[Chap.\ X, Theorem 7.4]{ML};
whence we may, and will, assume that $k$ contains all roots of unity.
In fact, then our computations show that the gap of zero Hochschild
cohomology groups from degree $3$ to the first non-zero group
thereafter can be arbitrarily large by choosing $q$ to be a primitive
$r$-th root of unity for large $r$.

Let $E(\L)$ be the Koszul dual of $\L$; that is,
$E(\L)=\Ext^*_\L(k,k)$. Our computations of $\HH^*(\L_q)$ lead to a
description of the graded centre of $E(\L_q)$, which is given in
section \ref{section:1.3}. 

We note that the Hochschild homology of $\L_q$ for $q$ in $k$, is
computed in \cite{Han}. It is shown there that $\dim_k\HH_n(\L_q)\geq
2$ for all $n\geq 0$. 

Throughout we consider right modules unless otherwise explicitly
stated. The enveloping algebra $\L^e$ of $\L$ is given by
$\L^\op\otimes_k \L$. 

\section{Background}

In this section we try to put Happel's question into a broader
context. We consider six conditions on a finite dimensional algebra
and describe how they are interrelated.

Let $\L$ be a finite dimensional algebra over a field $k$ with
Jacobson radial $\rrad$. The conditions we consider are the following.
\begin{enumerate}
\item[(1)] $\HH^n(\L)$ vanishes for $n$ sufficiently large.
\item[($1'$)] $\HH_n(\L)$ vanishes for $n$ sufficiently large.
\item[(2)] The global dimension of $\L$ is finite.
\item[(3)] For every field extension $k'$ of $k$, the global dimension
  of $k'\otimes_k\L$ is finite. 
\item[(4)] The projective dimension of $\L$ as a $\L^e$-module is
  finite. 
\item[(5)] The global dimension of $\L^e$ is finite.
\end{enumerate}
Happel's question is: Does (1) imply (2)? Our results show that this
is not true in general, even though it is true in the commutative case
\cite{AI}. There are some obvious connections; namely, that (3)
implies (2), (4) implies (1) and ($1^\prime$), and (5) implies (4). In
general, the following implications hold:
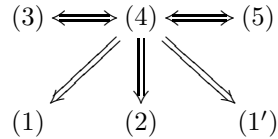
\[\xymatrix{
(3) \ar@{<=>}[r] & (4) 
\ar@{=>}[dl]\ar@{=>}[d]\ar@{=>}[dr]\ar@{<=>}[r] & (5)\\ 
(1) & (2) & (1')
}\]

The implication (4) implies (3) is a consequence of \cite[Corollary
IX.7.2]{CE}. If $k$ is algebraically closed, then Happel shows that 
(3) is equivalent to (4) \cite[1.5]{H2}. To see that (3) implies (4)
in general, consider a minimal projective resolution
$(\mathbb{P},\delta)$ of $\L$ over $\L^e$. If $k'$ is a field
extension of $k$, then $(k'\otimes_k\mathbb{P},1\otimes\delta)$ is a
minimal projective $(k'\otimes_k\L)^e$-resolution of $k'\otimes_k \L$,
since $k'\otimes_k \rrad\subseteq \rad(k'\otimes_k\L)$, and for a
projective $\L^e$-module $P$, the module $k'\otimes_k P$ is a
$(k'\otimes_k\L)^e$-module. Taking $k'$ to be the algebraic closure of
$k$, we see that if the projective dimension of $\L$ over $\L^e$ is
infinite, then the projective dimension of $k'\otimes_k\L$ over
$(k'\otimes_k\L)^e$ is also infinite. Applying Happel's result we
infer that not (4) implies not (3). 

The implication (4) implies (5) is found in \cite{vdB}. In that paper
van den Bergh suggests that condition (4) is a reasonable
non-commutative version of smoothness.

When $\L/\rrad$ is a separable algebra over $k$, the conditions (2),
(3), (4) and (5) are all equivalent. We have seen that (3), (4), and
(5) are equivalent and (3) implies (2) in general. To see that (2)
implies (4), we use the following two observations. First note that if
$\L/\rrad$ is a separable algebra over $k$, then $\L^e/\rad\L^e$ is
isomorphic to $\L/\rrad\otimes_k\L/\rrad\simeq
\Hom_k(D(\L/\rrad),\L/\rrad)$, where $D$ denotes as usual the
$k$-dual. Secondly, $\Ext^i_{\L^e}(\L,\Hom_k(D(\L/\rrad),\L/\rrad))$
is isomorphic to $\Ext^i_\L(D(\L/\rrad),\L/\rrad)$ by
\cite[Proposition IX.4.3]{CE}. In particular, if $\L/\rrad$ is a
separable algebra over $k$, (2) implies (1). In general (2) does not
imply (1) by the following example. Let $k$ be a field of
characteristic $2$, and let $\L=k(a)$ for some $a$ not in $k$ but
$a^2$ in $k$. By \cite{Ho} the Hochschild cohomology ring of $\L$ is
isomorphic to $k[x,y,z]/(x^2,y^2)$. Hence, the global dimension of
$\L$ is zero, but (1) fails for $\L$.

If $\L/\rrad$ is not separable over $k$, condition (2) does not imply
conditions (3), (4) and (5) in general. For example, if $\L$ is a
purely inseparable finite field extension of $k$ of the form $k(a)$,
then $\L^e$ is a local non-semisimple selfinjective algebra over $k$.
Hence the global dimension of $\L\otimes_k k(a)$ is infinite.

The results of this paper and \cite{Han} show that (1) does not imply
($1'$).  Han \cite{Han} also proves that ($1'$) is equivalent to (2) for
monomial algebras and for commutative algebras the equivalence was
shown by Avramov and Vigue\'e-Poirrier \cite{AV}. Han conjectures that
($1'$) and (2) are equivalent.

\section{The proof}

Let $\L=\L_q=k\langle x,y\rangle/(x^2,xy+qyx,y^2)$ for some field $k$
and for some element $q$ in $k\setminus\{0\}$. Order the elements by
length left lexicographic order by choosing $x>y>1$. The set
$\{x^2,xy+qyx,y^2\}$ is then a (non-commutative) Gr\"obner basis for
the ideal $(x^2,xy+qyx,y^2)$, and the normal forms of the elements in
$\L$ are given by $\{1,x,y,yx\}$ (see \cite{G}). Since $\L$ has a
quadratic Gr\"obner basis, $\L$ is a Koszul algebra (see \cite{GH}).

First we want to find a minimal projective resolution
$(\mathbb{P},\delta)$ of $\L$ over $\L^e$.  For each $n\geq 0$, we now
construct certain elements $\{ f^n_i\}_{i=0}^n$ in $\L^{\otimes_k n}$,
where $\L^{\otimes_k 0}$ means $k$. Let $f^0_0=1$, $f^1_0=x$ and
$f^1_1=y$. By letting $f^n_{-1}=0=f^n_{n+1}$ for all $n\geq 0$, we
define $\{f^n_i\}$ for $n\geq 2$ inductively by setting
\[f^n_i = f^{n-1}_{i-1}\otimes y+q^if^{n-1}_i\otimes x\]
for $i$ in $\{0,1,\ldots,n\}$.  It follows that $f^n_i$ is a linear
combination of all tensors of $x$'s and $y$'s of length $n$ with $y$
occurring $i$ times. The coefficient $q^\alpha$ in front of such a
$n$-fold tensor $w$ is determined by setting $\alpha$ equal to the
number of $x$'s that the $y$'s have to move across starting with the
$n$-fold tensor $\underbrace{x\otimes\cdots\otimes
x}_i\otimes\underbrace{y\otimes\cdots\otimes y}_{n-i}$ to obtain the
given tensor $w$. Using this we infer that
\[f^n_i = x\otimes f^{n-1}_i+q^{n-i}y\otimes f^{n-1}_{i-1}.\]

Let $P^n=\amalg_{i=0}^n \L\otimes_k f^n_i\otimes_k \L\subseteq
\L^{\otimes_k(n+2)}$ for $n\geq 0$, and let
$\tilde{f}^n_i=1\otimes f^n_i\otimes 1$ for $n\geq 1$ and
$\tilde{f}^0_0=1\otimes 1$. Define $\delta^n\colon
P^n\to P^{n-1}$ by letting 
\[\delta^n(\tilde{f}^n_i)=x\tilde{f}^{n-1}_i +
q^{n-i}y\tilde{f}^{n-1}_{i-1} + (-1)^n\tilde{f}^{n-1}_{i-1}y +
(-1)^nq^i\tilde{f}^{n-1}_ix.\]
Direct computations show that $\delta^2=0$.  

One can verify directly that $(\mathbb{P},\delta)$ is a projective
$\L^e$-resolution of $\L$ as follows. Let $\G$ be a Koszul $k$-algebra
where the degree zero part is a $k$. Suppose that $M$ is a Koszul
module over $\G$ for which the Betti numbers $\{ b_n\}_{n\geq 0}$ in a
minimal projective resolution are known, and that
\[\cdots\to G^n\extto{d^n} G^{n-1}\to\cdots \to G^1\extto{d^1}
G^0\extto{d^0} M\to 0\] is a linear complex $\mathbb{G}$ of projective
$\G$-modules $G^i$ generated by $b_i$ generators in degree $i$, where
$d^0\colon G^0\to M$ is onto and $d^i|_{G^i_i}\colon G^i_i\to
G^{i-1}_i$ is injective for all $i\geq 0$. Then it follows easily by
induction that $\mathbb{G}$ is a minimal projective resolution of $M$.

In our situation we apply this in the following way: Let $d^0\colon
P^0\to \L$ be the multiplication map. The algebra $\L^e$ is a Koszul
algebra over $k$. By a classical formula of Cartan-Eilenberg
$\Ext^n_{\L^e}(\L,k)\simeq \Ext^n_\L(k,k)$. The Koszul dual
$E(\L)=\Ext^*_\L(k,k)$ of $\L$ is isomorphic to $k\langle
x,y\rangle/(yx-qxy)$, so that the Betti numbers of a minimal
projective resolution of $\L$ over $\L^e$ are $\{ b_n=n+1\}_{n\geq
0}$. Moreover by \cite{GZ} (see \cite{BGSS} for a proof) a graded
algebra $\G$ generated in degrees $0$ and $1$ is Koszul if and only if
$\G$ is Koszul as a right module over $\G^e$. Hence $\L$ is a Koszul
module over $\L^e$. It is straightforward to check that the above
conditions are satisfied for $(\mathbb{P},\delta)$. Hence
$(\mathbb{P},\delta)$ is a minimal projective resolution of $\L$ over
$\L^e$.

Our next immediate goal is to obtain a multiplication formula for
elements in the Hochschild cohomology ring; see the formula \eqref{eq1}
below.  Recall that the bar resolution $\mathbb{B}$ of $\L$ over
$\L^e$ is given by $B^n=\L^{\otimes_k(n+2)}$ and differential
$d_n\colon B^n\to B^{n-1}$ given by
\begin{multline}
d_n(\lambda_0\otimes\lambda_1\otimes\cdots\otimes \lambda_{n+1}) \notag\\
=\sum_{i=0}^n
(-1)^i(\lambda_0\otimes\lambda_1\otimes\cdots\otimes\lambda_{i-1}
\otimes \lambda_i\lambda_{i+1}\otimes \lambda_{i+2}\otimes\cdots
\otimes \lambda_{n+1})\notag
\end{multline}
Define $\partial_i\colon B^n\to B^{n-1}$ by 
\[\partial_i(\lambda_0\otimes\cdots\otimes\lambda_{n+1})=(\lambda_0\otimes\cdots\otimes\lambda_{i-1}\otimes
\lambda_i\lambda_{i+1}\otimes\lambda_{i+2}\otimes\cdots\otimes\lambda_{n+1})\] for
$i=0,1,\ldots,n$. Then $d_n=\sum_{i=0}^n (-1)^n\partial_i$. 

We view $P^n$ as a submodule of $B^n$ in the natural way. For any
given $i$, consider the tensors of the form $w\otimes x\otimes
x\otimes w'$, $w\otimes x\otimes y\otimes w'$, $w\otimes y\otimes
x\otimes w'$, or $w\otimes y\otimes y\otimes w'$ occurring in a
$\tilde{f}^n_i$, where $w$ is an $t$-fold tensor. The first and the
last of these are mapped to zero by $\partial_t$.  If a tensor $q^tw\otimes
x\otimes y\otimes w'$ occurs in some $\tilde{f}^n_i$, then
$q^{t+1}w\otimes y\otimes x\otimes w'$ also occurs in
$\tilde{f}^n_i$. Then the sum of these elements is mapped to zero by
$\partial_t$. It follows from this that the element $\tilde{f}^n_i$ is mapped
to the element
\[x\tilde{f}^{n-1}_i+q^{n-i}y\tilde{f}^{n-1}_{i-1}+ 
(-1)^n\tilde{f}^{n-1}_{i-1}y + (-1)^nq^i\tilde{f}^{n-1}_ix\] 
by the differential in $\mathbb{B}$. This shows that
$(\mathbb{P},\delta)$ is a subcomplex of $(\mathbb{B},d)$.

The bar resolution $\mathbb{B}$ admits a comultiplication
$\Delta\colon \mathbb{B}\to \mathbb{B}\otimes_\L \mathbb{B}$, that is
a morphism of complexes over the identity on $\L$, given by
\[\Delta(\lambda_0\otimes\cdots\otimes\lambda_{n+1}) =
\sum_{i=0}^n
(\lambda_0\otimes\cdots\otimes\lambda_i\otimes1)\otimes_\L
(1\otimes\lambda_{i+1}\otimes\cdots\otimes \lambda_{n+1}).\] The (cup)
product in $\HH^*(\L)$ of two cycles $\eta$ and $\theta$ from
$\Hom_{\L^e}(\mathbb{B},\L)$ is given as the composition of the maps
\[\mathbb{B}\extto{\Delta} \mathbb{B}\otimes_\L\mathbb{B}
\extto{\eta\otimes \theta} \L\otimes_\L\L\extto{\nu}\L,\] where
$\nu\colon \L\otimes_\L\L\to \L$ is the natural map.  Let $\mu\colon
\mathbb{P}\to \mathbb{B}$ be the natural inclusion, and let $\pi\colon
\mathbb{B}\to \mathbb{P}$ be a chain map such that
$\pi\mu=\id_{\mathbb{P}}$.

Suppose that the image $\Delta(\mu\mathbb{P})$ is already contained in
$\mu\mathbb{P}\otimes_\L\mu\mathbb{P}\subseteq
\mathbb{B}\otimes_\L\mathbb{B}$. In that case, $\Delta$ induces a
comultiplication on $\mathbb{P}$ that we denote by $\Delta'$. In
particular $\Delta\mu=(\mu\otimes\mu)\Delta'$.  For homogeneous
elements $\eta$ and $\theta$ in $\HH^*(\L)$ we can represent these
elements as morphisms $\eta\colon P^m\to \L$ and $\theta\colon P^n\to
\L$. Using the bar resolution they can be represented by $\eta\pi_m$
and $\theta\pi_n$. Then we have
\begin{align}
\eta\ast \theta=(\eta\ast\theta )\pi_{m+n}\mu_{m+n} & = (\eta\pi_m\cup
			\theta\pi_n)\mu_{m+n}\notag\\ & =
			\nu(\eta\pi_m\otimes
			\theta\pi_n)\Delta\mu_{m+n}\notag\\ & =
			\nu(\eta\pi_m\otimes
			\theta\pi_n)(\mu_m\otimes\mu_n)\Delta'\notag\\ & =
			\nu(\eta\otimes \theta)\Delta'\notag
\end{align}
where $\eta\ast \theta$ denotes the product of $\eta$ and $\theta$ in
$\HH^*(\L)$.

Now we show that indeed $\Delta$ maps $\mu(\mathbb{P})$ to
$\mu\mathbb{P}\otimes_\L\mu\mathbb{P}$. Using the description of the
terms in $f^n_i$ it is straightforward to verify that
\[f^n_i=\sum_{j=\max\{0,i+t-n\}}^{\min\{t,i\}} q^{j(n-i+j-t)} 
f^t_j\otimes f^{n-t}_{i-j}\]
for $t=0,1,\ldots,n$. From this we infer that
\[\Delta(\tilde{f}^n_i) =
\sum_{t=0}^n\sum_{j=\max\{0,i+t-n\}}^{\min\{t,i\}}q^{j(n-i+j-t)}
\tilde{f}^t_j\otimes_\L \tilde{f}^{n-t}_{i-j},\] hence we see that
$\Delta\colon \mu\mathbb{P}\to \mu\mathbb{P}\otimes_\L \mu\mathbb{P}$
and we can define $\Delta'$ by the above formula. 

This enables us to describe the multiplication in $\HH^*(\L)$, which
we do next.  Let $\eta\colon P^m\to \L$ and $\theta\colon P^n\to \L$
represent elements in $\HH^m(\L)$ and $\HH^n(\L)$, respectively, where
$\eta(\tilde{f}^m_i)=\lambda_i$ for $i=0,1,\ldots,m$ and
$\theta(\tilde{f}^n_i)=\lambda_i'$ for $i=0,1,\ldots,n$. Then
\begin{equation}
\eta\ast \theta(\tilde{f}^{m+n}_i)=\sum_{j=\max\{0,i-n\}}^{\min\{m,i\}} 
q^{j(n-i+j)}\lambda_j\lambda_{i-j}'\label{eq1}
\end{equation}
for $i=0,1,\ldots,m+n$. We note that the product $\eta\ast\theta$ in
formula \eqref{eq1} is written as a map $P^{m+n}\to \L$, and hence may
be zero as an element in $\HH^*(\L)$.  In \cite{BGSS} we show that this
multiplicative structure can be generalized to any Koszul algebra
$\L=kQ/I$ for a finite quiver $Q$.\medskip

Now let us consider the vector space structure of $\HH^n(\L)$ for all
$n\geq 0$. 

Since $P^i\simeq (\L^e)^{i+1}$, we see that $\Hom_{\L^e}(P^{n-1},\L)$
is isomorphic to $\L^n$, which we view as an identification.  Let
$\eta\colon P^{n-1}\to \L$ be given by $(\lambda_0,\ldots,
\lambda_{n-1})$. Then for $n\geq 1$
\begin{multline}
\eta\delta^n  = (x\lambda_0+(-1)^n\lambda_0x,\ldots,\\
\underbrace{x\lambda_j + q^{n-j}y\lambda_{j-1}+(-1)^n\lambda_{j-1}y + 
(-1)^nq^j\lambda_jx}_{\text{$j$-th coordinate}},\ldots,\\
y\lambda_{n-1}+(-1)^n\lambda_{n-1}y)\notag
\end{multline}
The image $\Im(\delta^n)^*$ is spanned by vectors obtained from the
above by choosing $\lambda_i$ to be in the set $\{1,x,y\}$, since the
value $yx$ only gives zero. When $\lambda_j$ is $x$ or $y$ and all the
other $\lambda_i$'s are zero, then the result is a standard basis
vector times $ayx$ for some $a$ in $k$. Hence this can give rise to at
most $n+1$ linearly independent elements. When $\lambda_j=1$ and all
the other $\lambda_i$'s are zero, then the resulting element has the
form $(0,\ldots,0,\underbrace{ax}_j,by,0,\ldots,0)$ for some $a$ and
$b$ in $k$ and for $j=0,1,\ldots,n-1$. The non-zero elements of this
form are linearly independent and also linearly independent of those
obtained by choosing $\lambda_j$ in $\{x,y\}$. This shows that
$\dim\Im(\delta^n)^*\leq 2n+1$ for $n\geq 2$.

When $\lambda_i=0$ for $i\neq j$ and $\lambda_j=y$, then 
\[\eta\delta^n=(0,\ldots,0,\underbrace{(-q+(-1)^nq^j)yx}_j,0,\ldots,0)\]
for $j=0,1,\ldots,n-1$. If $\lambda_{j-1}=x$ and all other $\lambda_i$
are zero, then 
\[\eta\delta^n=(0,\ldots,0,\underbrace{(q^{n-j}-(-1)^nq)yx}_j,0,\ldots,0)\]
for $j=1,2,\ldots,n$. The coefficients in the $j$-coordinate of these
elements are zero if (i) for $j=0$ we have $q=(-1)^n$, (ii) for
$j=1,2,\ldots,n-1$ we have $q^{j-1}=(-1)^n$ and $q^{n-j-1}=(-1)^n$,
and (iii) for $j=n$ we have $q=(-1)^n$. For each $j$ such that one
of these three cases occurs, the dimension of $\Im(\delta^n)^*$
``decreases'' by one.  If $q$ is not a root of unity and $n\geq 3$,
then these conditions never hold.

When $\lambda_i=0$ for $i\neq j$ and $\lambda_j=1$, then 
\[\eta\delta^n=(0,\ldots,0,\underbrace{(1+(-1)^nq^j)x}_j,(q^{n-j-1}+(-1)^n)y,
0, \ldots,0)\] for $j=0,1,\ldots,n-1$. The coefficients in this vector
are zero if and only if $q^j=(-1)^{n+1}$ and
$q^{n-j-1}=(-1)^{n+1}$. If $q$ is not a root of unity and $n\geq 2$,
these equations are never satisfied. Hence, when $q$ is not a root
of unity and $n\geq 3$, the dimension of $\Im(\delta^n)^*$ is $2n+1$.

Using these facts, let us first compute the groups $\HH^n(\L)$ for
$n=0,1,2$. We have that
\[\HH^0(\L)=Z(\L)=
\begin{cases}
\Span_k\{1,yx\}, & q\neq -1\\
\L,        & q=-1
\end{cases}
\]
\sloppy From the above computations $\HH^1(\L)$ is spanned as a vector space
in $\Ker(\delta^2)^*/\Im(\delta^1)^*$ by the following maps 
\[
\begin{cases}
\{(x,0),(0,y)\}, & q\neq -1,1\\
\{(x,0),(yx,0),(0,yx),(0,y)\}, & q=-1, q\neq 1\\
\{(x,0),(y,0),(0,x),(0,y)\}, & q\neq -1, q=1\\
\Hom_{\L^e}(P^1,\L), & q=1, \kar k=2
\end{cases}\]
Similarly as above $\HH^2(\L)$ is spanned by the following maps 
\[
\begin{cases}
\{(0,yx,0)\}, & q\neq -1,1\\
\{(1,0,0),(y,0,0),(0,yx,0),(0,0,x),(0,0,1)\}, & q=-1, q\neq 1\\ 
\{(1,0,0),(0,1,0),(0,0,1),(yx,0,0),(0,yx,0),(0,0,yx)\}, & q\neq -1, q=1\\
\Hom_{\L^e}(P^2,\L), & q=1, \kar k=2
\end{cases}\]

\subsection{When $q$ is not a root of unity and non-zero} 
Suppose now that $q$ is not a root of unity. Since $\dim_k
(P^n)^*=4(n+1)$ and
$\dim_k\Im(\delta^n)^*+\dim_k\Im(\delta^{n+1})^*=4(n+1)$ for $n\geq
3$, it follows that $\Im(\delta^n)^*=\Ker(\delta^{n+1})^*$ and
$\HH^n(\L)=(0)$ for those values of $n$. By the considerations above
we have that (i) $\HH^0(\L)=\Span_k\{1,yx\}$, (ii)
$\HH^1(\L)=\Span_k\{(x,0),(0,y)\}$ and (iii)
$\HH^2(\L)=\Span_k\{(0,yx,0)\}$. Denote the element $yx$ in
$\HH^0(\L)$ by $z$ and the two elements $(x, 0)$ and $(0,y)$ in
$\HH^1(\L)$ by $u_0$ and $u_1$, respectively. Using the multiplication
worked out earlier, it is straightforward to compute that the elements
$\{z^2,zu_0,u_0z,zu_1,u_1z,u_0^2,u_1^2\}$ are all zero, while
$u_0u_1=-q(0,yx,0)$ and $u_1u_0=q(0,yx,0)$.  Hence $\HH^*(\L)$ is
isomorphic to the fibre product \[k[z]/(z^2)\times_k
\wedge^*(u_0,u_1),\] where $\wedge^*(u_0,u_1)$ is the exterior algebra
on the generators $u_0$ and $u_1$.

\subsection{When $q$ is a root of unity}
Suppose now that $q$ is a primitive $r$-th root of unity. Fix an
$n$. Whenever $\dim_k\Im(\delta^n)^*<2n+1$ the group $\HH^{i}(\L)\neq
(0)$ for $i$ in $\{n-1,n\}$. We observed above that this happened if 
\begin{enumerate}
\item[(a)]
\begin{enumerate}
\item[(i)] for $j=0$ we have $q=(-1)^n$, 
\item[(ii)] for some $j=1,2,\ldots,n-1$ we have $q^{j-1}=(-1)^n$ and
$q^{n-j-1}=(-1)^n$, 
\item[(iii)] for $j=n$ we have $q=(-1)^n$,
\end{enumerate}
\item[(b)] for some $j=0,1,\ldots,n-1$ we have $q^j=(-1)^{n+1}$ and
$q^{n-j-1}=(-1)^{n+1}$.
\end{enumerate}
When $j=1$ and $n-2=mr$ with $n$ even for some $m\geq 1$, then
condition (a)(ii) is satisfied. When $j=0$ and $n-1=mr$ with $n$ odd for
some $m\geq 1$, then condition (b) is satisfied. This shows that when
$q$ is a root of unity, then $\HH^*(\L)$ is infinite dimensional over
$k$.

The Koszul dual $E(\L)$ of $\L$ is isomorphic to $k\langle
x,y\rangle/(yx-qxy)$. Recall that the graded centre $Z_\gr(E(\L))$ of
$E(\L)$ is generated by $\{z\in E(\L)\mid z \text{\ homogeneous,\ }
z\gamma=(-1)^{\deg z\deg\gamma}\gamma z \text{\ for all homogeneous\ }
\gamma\in E(\L)\}$. Using that $E(\L)$ is bigraded it is
straightforward to compute $Z_\gr(E(\L))$. If $r$ is even or $\kar
k=2$, then $Z_\gr(E(\L))$ is equal to $H'=k[x^r,y^r]$. Otherwise
$Z_\gr(E(\L))$ is equal to $k[x^{2r},x^ry^r,y^{2r}]$ and therefore
contains $H'=k[x^{2r},y^{2r}]$. In either case $E(\L)$ is a finitely
generated module over $H'$. By \cite{BGSS} or \cite{K} the image of the
natural map $\HH^*(\L)\to E(\L)$ is $Z_\gr(E(\L))$. Take homogeneous
inverse images of the two generators of $H'$, and let $H$ denote the
polynomial subring of $\HH^*(\L)$ they generate. Since $\E(\L)\simeq
\HH^*(\L,S)$ where $S$ is the only simple $\L^e$-module, $\HH^*(\L)$
is also a finitely generated module over $H$ as $H$ is noetherian and
$\L$ is filtered in $S$ as a $\L^e$-module.

\subsection{The graded centres of $E(\L_q)$}\label{section:1.3}

By direct calculations and the remarks above we have the following
description of the graded centre of $E(\L_q)$ for all $q$ in $k$. 

\begin{prop}
The graded centre of $E(\L_q)$ is as follows.
\[Z_\gr(E(\L_q))=\begin{cases}
k, & q \text{\ is not a root of unity}\\
k[x^r,y^r], & \begin{cases} 
q \text{\ is a primitive $r$-th unity and $r$ even, or}\\
\kar k=2
\end{cases}\\
k[x^{2r},x^ry^r,y^{2r}], & \text{otherwise}
\end{cases} \]
\end{prop}

\section{Further details}
In this section we give the full ring structure of $\HH^*(\L_q)$ and
the Hilbert series with respect to the cohomological degree, when $q$
is an $r$-th primitive root of unity and for completeness when $q=0$.

The description of the ring structure of $\HH^*(\L)$ is divided into
several cases. The vector space basis for $\HH^i(\L)$ for $i=0,1,2$
are given in the previous section when $q\neq 0$ and below when
$q=0$. We first give a vector space basis for $\HH^i(\L)$ for $i\geq
3$ in the different cases, and then the full ring structure of
$\HH^*(\L)$. 

The cases we have to consider are the following (a) $r>1$ is odd and
$\kar k\neq 2$, (b) $\kar k=2$ with $q\neq 1$, or $\kar k\neq 2$ and
$r>2$ is even, (c) $\kar k=2$ and $q=1$, (d) $\kar k\neq 2$ and $q=-1$
and finally (e) $\kar k\neq 2$ and $q=1$.

\subsection{The case $r>1$ is odd and $\kar k\neq
2$}\label{section:2.1} 
In this case representing maps for $\HH^n(\L)$ when $n>2$ are given by
\[\HH^n(\L)=
\begin{cases}
\Span_k\{ e_j\}_{j\in\{tr\}_{t=0}^{2s}}, & n=2sr\\ \Span_k\{\{
xe_{j-1}\}_{j\in\{tr+1\}_{t=0}^{2s}}, \{
ye_j\}_{j\in\{tr+1\}_{t=0}^{2s}}\}, & n=2sr+1\\ \Span_k\{
yxe_j\}_{j\in\{tr+1\}_{t=0}^{2s}}, & n=2sr+2\\ 0, & \text{otherwise}
\end{cases}\]
Using the description of the multiplication given in the previous
section and remembering that $q$ is an $r$-th root of unity, one can 
see that 
\[\HH^*(\L)\simeq k[z]/(z^2)\times_k 
\left(\wedge^*(u_0,u_1)[w_0,w_1,w_2]/(w_0w_2-w_1^2)\right)\] where $z$
is in degree $0$, $u_0=(x,0)$ and $u_1=(0,y)$ are in degree $1$, and
$w_i\colon P^{2r}\to \L$ is in degree $2r$ and given by
$w_i(e_{ir})=1$ and zero otherwise for $i=0,1,2$. We have that
$\sum_{i\geq 0} \dim_k \HH^i(\L) t^i =
1+\frac{(1+t)^2(1+t^{2r})}{(1-t^{2r})^2}$.

\subsection{The case $\kar k=2$ with $q\neq 1$, or $\kar k\neq 2$ and
$r>2$ is even} Representing maps for $\HH^n(\L)$ when $n>2$ are in
this case given by
\[\HH^n(\L)=
\begin{cases}
\Span_k\{ e_j\}_{j\in\{tr\}_{t=0}^s}, & n=sr\\
\Span_k\{\{ xe_{j-1}\}_{j\in\{tr+1\}_{t=0}^s},
\{ ye_j\}_{j\in\{tr+1\}_{t=0}^s}\}, & n=sr+1\\
\Span_k\{ yxe_j\}_{j\in\{tr+1\}_{t=0}^s}, & n=sr+2\\
0, & \text{otherwise}
\end{cases}\]
Similar as above we find that 
\[\HH^*(\L)\simeq k[z]/(z^2)\times_k
\left(\wedge^*(u_1,u_2)[w_0,w_1]\right)\] where $z$ is in degree $0$,
$u_0=(x,0)$ and $u_1=(0,y)$ are in degree $1$, and $w_i\colon P^r\to
\L$ is in degree $r$ and given by $w_i(e_j)=1$ for $j=ir$ and zero
otherwise for $i=0,1$. The Hilbert series is given by $\sum_{i\geq 0}
\dim_k \HH^i(\L) t^i = 1 +\frac{(1+t)^2}{(1-t^r)^2}$.  

\subsection{The case $\kar k=2$ and $q=1$} Here the differential in
the complex $\Hom_{\L^e}(\mathbb{P},\L)$ is zero, so that all maps in
$\Hom_{\L^e}(P^n,\L)$ represent non-zero elements in $\HH^*(\L)$.  We
have that $\HH^*(\L)\simeq \L[w_0,w_1]$ where $w_i\colon P^1\to \L$ is
given by $w_i(e_j)=1$ for $i=j$ and zero otherwise for
$i=0,1$. Furthermore $\sum_{i\geq 0} \dim_k \HH^i(\L) t^i =
\frac{4}{(1-t)^2}$. 

\subsection{The case $\kar k\neq 2$ and $q=-1$}
In this case representing maps for $\HH^n(\L)$ when $n>2$ are given by
\[\HH^n(\L)=
\begin{cases}
\Span_k\{\{ e_j\}_{j\in\{2t\}_{t=0}^{\frac{n}{2}}}, ye_0, xe_n, 
\{ yxe_j\}_{j\in\{2t-1\}_{t=1}^{\frac{n}{2}}}\}, & n \text{\ even}\\
\Span_k\{\{ xe_{j-1}\}_{j\in\{2t-1\}_{t=1}^{\frac{n+1}{2}}}, 
\{ ye_j\}_{j\in\{2t-1\}_{t=1}^{\frac{n+1}{2}}}, 
yxe_0, yxe_n\}, & n \text{\ odd}
\end{cases}\]
This gives 
\[\HH^*(\L)\simeq 
\left(\L\otimes_k
\wedge^*(u_0,u_1)\right)[w_0,w_1]/(xu_0,yu_1,xw_0,yw_1),\] 
where $u_0=(x,0)$ and $u_1=(0,y)$ are in degree $1$ and $w_0=(1,0,0)$
and $w_1=(0,0,1)$ are in degree $2$. The Hilbert series is given by
$\sum_{i\geq 0} \dim_k \HH^i(\L) t^i = \frac{4-4t+t^2}{(1-t)^2}$.

As $\L_{-1} \simeq k[x]/(x^2)\otimes_k k[y]/(y^2)$, there is an
isomorphism of vector spaces 
\begin{align}
\HH^*(\L_{-1}) & \simeq \HH^*(k[x]/(x^2))\otimes_k \HH^*(k[y]/(y^2))\notag\\
& \simeq (k[x]/(x^2) \times_k \wedge^*(u_0)[w_0])\otimes_k (k[y]/(y^2)
\times_k \wedge^*(u_1)[w_1])\notag
\end{align}
by \cite[Chap.\ X, Theorem 7.4]{ML}, which is in accordance with our result!

\subsection{The case $\kar k\neq 2$ and $q=1$}\label{section:2.5}
In this case representing maps for $\HH^n(\L)$ when $n>2$ are given by
\[\HH^n(\L)=
\begin{cases}
\Span_k\{\{ e_j\}_{j=0}^n, \{ yxe_j\}_{j=0}^n\}, & n \text{\ even}\\
\Span_k\{\{ xe_j\}_{j=0}^n, \{ ye_j\}_{j=0}^n\}, & n \text{\ odd}
\end{cases}\]
It follows that
\[\HH^*(\L)\simeq \Big(k[z]/(z^2)\times_k
\wedge^*(u_0,u_1,u_2,u_3)\Big)[w_0,w_1,w_2]/I,\] where $z$ is in
degree $0$, $u_0=(x,0)$, $u_1=(y,0)$, $u_2=(0,x)$, $u_3=(0,y)$ are in
degree $1$, and $w_i\colon P^2\to \L$ is given by $w_i(e_j)=1$ for
$j=i$ and zero otherwise for $i=0,1,2$. The ideal $I$ is given by
\begin{multline}
(u_0u_2,u_1u_3, u_0u_1+zw_0, u_0u_3+zw_1, u_2u_3+zw_2, u_1u_2-zw_1,\notag\\ 
u_0w_1-u_2w_0, u_1w_1-u_3w_0, u_0w_2 - u_2w_1, u_1w_2-u_3w_1,w_0w_2-w_1^2) 
\end{multline}
We have that $\sum_{i\geq 0} \dim_k \HH^i(\L) t^i = \frac{2}{(1-t)^2}$.

\subsection{The case $q=0$}\label{section:q=0} For completeness we now
consider the degenerate case when $q=0$. We have that the centre of
$\L$ is given by $\Span_k\{1,yx\}$, so that
$\HH^0(\L)=\Span_k\{1,yx\}$. It is left to the reader to verify that
$\HH^1(\L)=\Span_k\{(x,0),(0,y)\}$. For $n>1$ the representing maps
for $\HH^n(\L)$ are
\[\Span_k\{\{ xe_i\}_{i=0}^{n-2}, ye_n, \{ yxe_i\}_{i=1}^{n-1}\}.\]
It follows from this that all products in $\HH^{\geq 1}(\L)$ are
zero. Hence the Hochschild cohomology ring of $\L$ is isomorphic to
\[k[v_{00},v_{10},v_{11},\{v_{ij}\}_{i\geq 2,
j\in\{0,1,\ldots,2i-2\}}]/(v_{ij}v_{i'j'}),\] where $v_{ij}$ is in
degree $i$. Finally $\sum_{i\geq 0} \dim_k \HH^i(\L) t^i =
\frac{1+t^3}{(1-t)^2}$.

\end{document}